\documentclass[11pt,letter]{amsart}
\usepackage[utf8]{inputenc}
\usepackage[english]{babel}
\usepackage[dvipsnames]{xcolor}
\usepackage{rotating}
\usepackage[all]{xy}
\usepackage{comment}
\usepackage{pdflscape}
\usepackage{hyperref}
\usepackage{tikz}
\usepackage{tikz-cd}
\usepackage{soul,color}
\usepackage[colorinlistoftodos]{todonotes}
\usepackage[width=\textwidth,font=footnotesize,labelfont=bf, aboveskip=30pt]{caption}
\usepackage[symbol]{footmisc}
\usepackage[OT2,T1]{fontenc}
\DeclareSymbolFont{cyrletters}{OT2}{wncyr}{m}{n}
\DeclareMathSymbol{\Sha}{\mathalpha}{cyrletters}{"58}

\newcommand*{\fullref}[1]{\hyperref[{#1}]{\ref*{#1}. \nameref*{#1}}}

\usepackage[osf]{mathpazo}
\usepackage{charter}
\usepackage[top=3.2cm, bottom=3.2cm,left=2.8cm, right=2.8cm, heightrounded,bindingoffset=0mm]{geometry}

\usepackage{hyperref}
\hypersetup{bookmarksnumbered=true, linkcolor=due, citecolor=Green, urlcolor=black,}

\makeatletter
\@namedef{subjclassname@2020}{%
  \textup{2020} Mathematics Subject Classification}
\makeatother

\usepackage{fancyhdr}

\usepackage{amsfonts}
\usepackage{amssymb}
\usepackage{extarrows}
\usepackage{cleveref}
\usepackage{amsthm}
\usepackage{amscd}
\usepackage{amsmath}
\usepackage{mathtools}
\usepackage{mathrsfs}
\usepackage{float}

\usepackage{color}	
\definecolor{due}{RGB}{0,76,147}

\usepackage{graphicx}	
\usepackage{multicol}	
\usepackage{wrapfig}
\usepackage[english]{babel} 
\usepackage[utf8]{inputenc}
\usepackage[T1]{fontenc}
\usepackage{enumitem}

\usepackage{longtable}
\usepackage{cite}

\theoremstyle{definition}
\newtheorem{defi}{Definition}[section]
\theoremstyle{plain}
\newtheorem{thm}[defi]{Theorem}

\newtheorem{introthm}{Theorem}[section]

\newtheorem{introcor}[introthm]{Corollary}
\newtheorem{prop}[defi]{Proposition}
\newtheorem{cor}[defi]{Corollary}
\newtheorem{lemma}[defi]{Lemma}
\theoremstyle{remark}

\theoremstyle{definition}

\newcommand{\Pic}{\operatorname{Pic}}

\newcommand{\codim}{\operatorname{codim}}

  \makeatletter
\newcommand{\xdashrightarrow}[2][]{\ext@arrow 0359\rightarrowfill@@{#1}{#2}}
\newcommand{\xdashleftarrow}[2][]{\ext@arrow 3095\leftarrowfill@@{#1}{#2}}
\newcommand{\xdashleftrightarrow}[2][]{\ext@arrow 3359\leftrightarrowfill@@{#1}{#2}}
\def\rightarrowfill@@{\arrowfill@@\relax\relbar\rightarrow}
\def\leftarrowfill@@{\arrowfill@@\leftarrow\relbar\relax}
\def\leftrightarrowfill@@{\arrowfill@@\leftarrow\relbar\rightarrow}
\def\arrowfill@@#1#2#3#4{%
  $\m@th\thickmuskip0mu\medmuskip\thickmuskip\thinmuskip\thickmuskip
   \relax#4#1
   \xleaders\hbox{$#4#2$}\hfill
   #3$%
}
\makeatother

\newcommand{\escorta}[3]{0\rightarrow{#1}\longrightarrow{#2}
\longrightarrow{#3}\rightarrow0}
\newcommand{\fas}[1]{\mathscr{#1}}

\newcommand{\F}{\mathcal{F}}

\newcommand{\OO}[1]{\fas{O}_{#1}}

\numberwithin{equation}{section}
\begin{document}
	\title[Syzygies of Projective Bundles]{Cohomological rank functions and syzygies of projective bundles on abelian varieties}
	\author{Minyoung Jeon}
	  \address{Department of Mathematics\\University of Georgia, Athens GA 30602, USA}
\email{\url{minyoung.jeon@uga.edu}}
  \author{Sofia Tirabassi}
  \address{Department of Mathematics\\ Stockholm University, Albano campus Hus 1, Stockholm, Sweden}
\email{\url{tirabassi@math.su.se}}

%

\thanks{All authors contributed equally to this work.}

\clearpage\maketitle 

\begin{abstract}
We use Jiang--Pareschi's cohomological rank functions and techniques developed by Caucci and Ito to study syzygies of projective bundles on abelian varieties. Consequently, we generalize and refine Chiantipalli's results regarding property $N_p$ of projective bundles over abelian varieties.

\end{abstract}

\section{Introduction}

Let $X$ be an abelian variety over the field of complex numbers $\mathbb{C}$. 
Given $(X,l)$ a polarized abelian variety and an object $\mathcal{F}\in D^b(X)$, Jiang--Pareschi \cite{JiPa2020} and Caucci \cite{Caucci} define the {cohomological rank functions} associated to $\mathcal{F}$ and $l$
\[
h^i_{\mathcal{F},l}:\mathbb{Q}\rightarrow\mathbb{Q}^+
\]
and use them to introduce certain numerical invariants associated to the polarization $l$. For example, if $\mathcal{\F}=\mathcal{I}_x$ is the ideal sheaf of a point $x\in X$, the \emph{basepoint-freeness threshold} is defined as
\[\beta(l):=\inf\{x\in\mathbb{{Q}}^+\:|\:h^1_{\mathcal{I}_x,l}(x)=0\},\] 
and has the following properties:
\begin{itemize}
    \item $\beta(l)\leq 1$ and $\beta(l)<1$ if and only if all line bundles $L$ with $c_1(L)=l$ are globally generated;
    \item by definition (see \ref{crf} and \cite[Proof of Cor. 1.2]{Caucci}), we have that, for every $m$ positive integer, $\beta(ml)=\frac{\beta(l)}{m}$;
    \item $\beta(l)<\frac{1}{2}$ if, and only if, for every $L$ with $c_1(L)=l$, we have that $L$ is normally generated (see \cite[Cor. E]{JiPa2020});
    \item if $\beta(l)<\frac{1}{p+2}$, then, for every $L$ with $c_1(L)=l$, we have that $L$ satisfies property $N_p$ (see \cite[Thm. 1.1]{Caucci});
\end{itemize}
We recall that a vector bundle $E$ on an abelian variety $X$ is called homogeneous if, and only if, $t_x^*E\simeq E$ for all $x\in A$, where $t_x:X\rightarrow X$ denotes the "translation by $x$". This condition is known to be equivalent to $E$ having a filtration

\[
0=E_0\subset E_1\subset \cdots \subset E_\mathfrak{r}=E
\]
of subvector bundles with $\mathrm{rk}(E_i)=i$  and
\[
E_t/E_{t-1}\in\mathrm{Pic}^0(X)\;\text{for}\;1\leq t\leq\mathfrak{r}.
\] 
  Our main result is the following.
\begin{introthm}\label{TheoremA}
Let $X$ be an abelian variety of dimension $g$ over the complex number $\mathbb{C}$ and consider an homogeneous vector bundle $E$ on $X$.
Let $Y:=\mathbb{P}_X(E)$ be the associated projective vector bundle with a projective morphism $\pi:Y\rightarrow X$. We denote by $\mathcal{O}_Y(1)$ the tautological line bundle on $Y$. Let $p\geq 0$. Let $A$ be an ample line bundle on $X$ with 
\[\beta(c_1(A))<\min\left\{\frac{r+1}{p+r+2},\frac{1}{2}\right\}.\] Then $\mathcal{O}_Y(1)\otimes\pi^* A$ is very ample and satisfies property $N_p^r$.

In particular, if $\beta(c_1(A))<1/(p+2)$, then $\mathcal{O}_Y(1)\otimes\pi^* A$ satisfies property $N_p$.

\end{introthm}

Property $N_p^0$ is equivalent to property $N_p$. So, our first main theorem, Theorem \ref{TheoremA} generalizes Chintapalli's main result on property $N_p$ \cite[Theorem 1.2]{Chi19}, however we were told by A. Ito that Theorem A above is a particular case of the main theorem in the appendix of \cite{Raychaudhury+2024+19+37}, where the filtration assumption is not needed. As a corollary (see Corollary \ref{introcor1}), we get Chiantipalli original result that $\mathcal{O}_Y(1)\otimes\pi^* A^{\otimes n}$ satisfies $N_p$ for every $n\geq p+3$.

Projective normality and normal presentation on curves were initially studied by Castelnuovo, and later by Mattuck, Fujita and St-Donat. Mark Green merged these concepts as special cases of the $N_p$-properties, $p$-th syzygy property for $p\geq 0$, see \cite{green1984koszul,GreenII,GreenIII}. In particular, Green generalized the result by Castelnuovo, showing that a line bundle of degree at least $2g+1+p$ on a smooth curve of genus $g$ satisfies $N_p$-property. In case of higher dimensional varieties, Ein and Lazarsfeld \cite{EL}, among many others, proved that adjoint linear series associated to a very ample line bundle on a smooth projective variety satisfies the property $N_p$. Park \cite{Park} generalized the result by Ein--Lazarsfeld to projective bundles over smooth projective varieties, under certain hypothesis.
  
As for abelian varieties, for an ample line bundle $A$ on an abelian variety, Lazarsfeld made a conjecture that $A^n$ satisfies $N_p$-properties if $n\geq p+3$ for $p\geq 0$. This conjecture was proved and generalized to property $N_p^r$ by Pareschi \cite{pareschi2000syzygies} in characteristic $0$, where it is shown that $A^n$ satisfies $N_p^r$-property if $n\geq (p+r+3)/(r+1)$. Chintapalli and Iyer \cite{ChiIyer} proved the Pareschi result on property $N_p$ for hyperelliptic varieties, and Chintapalli \cite[Theorem 1.2]{Chi19} extended the result to homogeneous projective bundles over abelian varieties, that is to projectivizations of homogenous vector bundles.

 The theory of {\it M-regularity} is introduced and developed by Pareschi and Popa \cite{PP03,PP04,PP11}. Jiang and Pareschi \cite{JiPa2020} extend this notion to a $\mathbb{Q}$-{\it twisted sheaf} $\mathcal{F}\langle x\ell\rangle$ for a coherent sheaf on an abelian variety $X$, where $x$ is a rational number in $\mathbb{Q}$. See \S\ref{sec2.2} for the definitions of M-regularity and $\mathbb{Q}$-twisted sheaves. 
  
  The following is our second main theorem.  

\begin{introthm}[Theorem \ref{thm:5.4}]\label{TheoremC}
Suppose $(Y,X,\pi)$ as before. 
Let $A$ be an ample line bundle on $X$. Let $p\geq 1$.
 If $\mathcal{I}_o\left\langle \frac{r+1}{p+r+2}c_1(A)\right\rangle$ is M-regular, then $\mathcal{O}_Y(1)\otimes\pi^*A$ satisfies $N_p^r$-property.  
 
 In particular, if $\mathcal{I}_o\left\langle \frac{1}{p+2}c_1(A)\right\rangle$ is M-regular, then $\mathcal{O}_Y(1)\otimes\pi^*A$ satisfies $N_p$ property. 
\end{introthm}

We utilize the work of Ito \cite{Ito22} on the M-regularity of $\mathbb{Q}$-twisted sheaves to prove Theorem \ref{TheoremC}.

   Ito \cite{Ito22} refines results of Pareschi--Popa \cite{PP03,PP11} on global generation and surjectivity of multiplication maps of global sections of coherent sheaves on abelian varieties, using the theory of M-regularity of $\mathbb{Q}$-twisted sheaves. Furthermore, he shows the implication of the M-regularity of a suitable $\mathbb{Q}$-twisted sheaf to property $N_p$ and jet-ampleness for an ample line bundle on abelian varieties.

As a corollary of Theorem \ref{TheoremC}, we have the following:
\begin{introcor}\label{introcor2}
  Let $X$ an abelian variety, $E$ a vector bundle as in the statement of Theorem \ref{TheoremC}, and $A$ an ample line bundle on $X$ with no base divisor. If $n\geq (p+r+2)/(r+1)$, then $\mathcal{O}_Y(1)\otimes\pi^* A^{\otimes n}$ satisfies $N_p^r$ for $r\geq 0$, $p\geq 1$.
  \end{introcor}

  Our Corollary \ref{introcor2} gives refined results of Chintapalli, since it implies that $\mathcal{O}_Y(1)\otimes\pi^*A^{\otimes n}$ satisfies $N_p$ for $n\geq p+2$.

Let $X$ be a nodal curve, with one node, and $X_0$ its normalization. Desingularizations $\widetilde{J}(X)$ of compactified Jacobian $\overline{J}(X)$ can be expressed as a $\mathbb{P}^1$-bundle $\mathbb{P}(E)$ over  the Jacobian $J(X_0)$, where $E$ is given by direct sum of degree zero line bundles on desingularization ${J}(X_0)$. In \cite[Theorem 5.4]{Chi19}, Chintapalli applied his result, Corollary \ref{introcor1} to desingularization of compactified Jacobians by setting $X=J(X_0)$ with the direct sum $E$ of degree zero line bundles.
 
  In the last section, \S\ref{sec5}, we provide an application of our result to desingularization of compactified Jacobians, generalizing \cite[Theorem 5.4]{Chi19}.

  This paper is organized as follows. In \S\ref{sec2}, we review some necessary definitions and notations. In \S\ref{sec3}, we recall ampleness and conditions for tautological line bundles to satisfy property $N_p^r$. In \S\ref{sec4}, we prove our maint results, Theorem \ref{TheoremA}, Theorem \ref{TheoremC}, and Corollary \ref{introcor2}. In \S\ref{sec5}, we study an application of our results to desingularization of compactified Jacobians.
  
\subsection*{Notation}
We work on the field of complex numbers $\mathbb{C}$. We use standard notation for abelian varieties. If $X$ is an abelian variety, we denote by $\hat{X}$ its dual abelian variety. If $X$ is an abelian variety, we denote by $\hat{X}$ its dual abelian variety.  Denote by $\mathscr{P}$ the reduced Poincar\'e bundle on $X\times\hat{X}$. If we view $\hat{X}$ as the moduli spaces of topologically trivial line bundles on $X$, we have that $\mathscr{P}$ is the unique universal family such that $\mathscr{P}_{\{0\}\times\hat{X}}\simeq\mathcal{O}_{\hat{X}}$ and $\mathscr{P}_{X\times\{\mathcal{O}_{X}\}}\simeq\mathcal{O}_{{X}}$. Given an element $\alpha\in\hat{X}$, we let $P_\alpha$ be the corresponding line bundle in $\operatorname{Pic}^0(X)$, that is 
$P_\alpha=\mathscr{P}_{|X\times\{\alpha\}}$. Throughout this paper, we may use $\alpha$ to indicate $P_\alpha$. Given a integer $n$, we let $n_X:X\rightarrow X$ be the "multiplication by $n$" map. When $n$ is positive, we denote by $X[n]\coloneqq \operatorname{Ker}n_X$: the group-subscheme of $n$-torsion points of $X$. For every $x\in X$ we let $t_x:X\rightarrow X$ be the map defined by $a\mapsto a+x$. Given a polarization $l$ on $X$, we let $\varphi_l:X\rightarrow \hat{X}$ be the isogeny defined by $x\mapsto t_x^*L\otimes L^{-1}$, where $L$ is an ample line bundle with $c_1(L)=l$. It is a standard fact about abelian varieties that this does not depends from the choice of $L$. We will denote the  kernel of $\varphi_l$ by $K(l)$. The identity element of $X$ is denoted simply by 0.\\ By $D^b(X)$ we denote $D^b(\operatorname{Coh}(X))$, the derived category of the category of coherent sheaves on $X$.

\section{Preliminaries}\label{sec2}
\subsection{Syzygies of projective varieties}\label{sub:N_P}
Let $Z$ be an algebraic variety over $\mathbb{C}$ and let $L$ be an ample
invertible sheaf on $Z$. Let $S_L$ be the symmetric algebra over $H^0(Z,L)$. Then \textit{section ring associated to  $L$},
\[R_L\coloneqq \bigoplus_{m\in\mathbb{Z}}H^{0}({Z},{L^{
\otimes m}}),\]
is a finitely generated graded $S_L$-module and as such it admits a \textit{minimal free
resolution} 
\begin{equation}\label{minalfree}
F_\bullet = 0\rightarrow\cdots\xrightarrow{f_{p+1}}
F_p\xrightarrow{f_p}\cdots\xrightarrow{f_2} F_1\xrightarrow{f_1}
F_0\xrightarrow{f_0}R_L\rightarrow 0
\end{equation}
with  $F_i\simeq\bigoplus_j S_L(-a_{ij})$, $a_{ij}\in\mathbb  Z_{>0}$.
Following Green \cite{green1984koszul},  we say that $L$ satisfies
property $N_p$ for some nonnegative integer $p$, if, in the notations above,
$$F_0=S_L$$
and
$$F_i=\oplus S_L(-i-1)\quad 1\leq i\leq p.$$
In \cite{pareschi2000syzygies}, Pareschi  extended this by introducing property $N_p^r$, where $r\geq 0$ is an integer which measures how much property $N_p$ fails (that is property $N_p^0$ is equivalent to property $N_p$). More precisely, we say that, $L$ satisfies property $N_p^r$ if, in the notation above, $a_{0j}\leq 1+r$ for
every $j$. Inductively, we see that property $N_p^r$ holds for $L$ if $L$ satisfies  property $N_{p-1}^r$ and $a_{pj}\leq p+1+r$ for every $j$.\par
 It is well known that condition $N_p$ is equivalent to the exactness in the middle of the complex
\begin{equation}\label{eq:np1}
 \bigwedge^{p+1}H^0(L)\otimes H^0(L^{\otimes h})\rightarrow \bigwedge^{p}H^0(L)\otimes H^0(L^{\otimes
h+1})\rightarrow \bigwedge^{p-1}H^0(L)\otimes H^0(L^{\otimes h+2})
\end{equation}
for any $h\geq 1$. More  generally, condition $N_p^r$ is equivalent to exactness in the middle of
\eqref{eq:np1} for every $h\geq r+1$. When $L$ is globally generated, we can consider the following exact sequence:
\begin{equation}
 \escorta{M_L}{H^0(L)\otimes\OO{Z}}{L}.
\end{equation}
Caucci in \cite[Prop.4.1]{Caucci} has shown that, independently from the the characteristic of the base field, we have that property $N_p$ (respectively $N_p^r$) is implied by the vanishing 
\begin{equation}\label{importan vanishing}
    H^1(Z,M_{L}^{\otimes p+1}\otimes L^{\otimes h})=0
\end{equation}
for every $h\geq 1$ (respectively, for every $h\geq r+1$). Therefore, our problem is reduced to check the vanishing of certain cohomology groups.
\subsection{$\mathbb{Q}$-twisted sheaves}\label{sec2.2}
Given a polarization $l$ on an abelian variety $X$ we define \emph{coherent object $\mathbb{Q}$-twisted by $l$} as equivalence classes pairs $(\mathcal{F}, xl)$ where $\mathcal{F}$ is an object in $D^b(X)$ and $x$ is a rational number, with the equivalence relation generated by
$$
(\mathcal{F}\otimes L^m,xl)\sim(\mathcal{F},(m+x)l)
$$
for every integer $m$, and every line bundle $L$ with $c_1(L)=l$. The equivalence class of $(\mathcal{F},xl)$ is denoted by $\mathcal{F}\langle xl\rangle$. Given a $\mathbf{Q}$-twisted sheaf $\mathcal{F}\langle xl\rangle$ we choose a representative $L$ for $l$ and define its cohomological support loci as
$$
V^i\left(X, \mathcal{F}\left\langle \frac{a}{b}l\right\rangle, L\right)\coloneqq\{  \alpha\in\hat{X}\:|\:h^i(X, b_X^*\mathcal{F}\otimes L^{ab}\otimes P_  \alpha)\neq0\}.
$$

If we change the representative of $l$ we obtain a new locus that is a translate of $V^i\left(X, \mathcal{F}\langle \frac{a}{b}l\rangle, L\right)$ by an element in $\hat{X}$. In particular, the dimension of these loci does not depend from the choice of the line bundle $L$.
If we change the representation of $\frac{a}{b}$, we have that still the dimension of these loci is unchanged.
 In fact, let $k$ be a positive integer. We have that 
 \begin{align*}
  (k b)_X^*\mathcal{F}\otimes L^{k^2ab} &\simeq   k_X^*b_X^* \mathcal{F}\otimes\left( L^{ab}\right)^{k^2}\\
  &\simeq k_X^*\left(b_X^* \mathcal{F}\otimes L^{ab}\otimes\beta\right)
 \end{align*}
 for $\beta\in\hat{X}$ such that $(-1)_X^*L^{\frac{k^2-k}{2}}\otimes\beta^k \simeq L^{\frac{k^2-k}{2}}$. We deduce that $V^i\left(X, \mathcal{F}\left\langle \frac{ka}{kb}l\right\rangle, L\right)$ is the image of a translate of $V^i\left(X, \mathcal{F}\left\langle \frac{a}{b}l\right\rangle, L\right)$ through the isogeny $k_X^*$.\par

 Thus we can introduce the following definition.
\begin{defi}\label{def2.1}
 Let $x$ a rational number. We say that $\mathcal{F}\left\langle xl\right\rangle$ is GV if $\codim V^i(X,\mathcal{F}\left\langle \frac{a}{b}l\right\rangle, L)\ge i$ for every $i\ge 1$, and for every representation $x=\frac{a}{b}$.\par
We say that $\mathcal{F}\left\langle xl\right\rangle$ is M-regular if $\codim V^i(X,\mathcal{F}\left\langle \frac{a}{b}l\right\rangle, L)\ge i+1$ for every $i\ge 1$, and for every representation $x=\frac{a}{b}$.\par
We say that a $\mathbb{Q}$-twisted sheaf $\mathcal{F}\left\langle \frac{a}{b}l\right\rangle$ satisfies I.T.(0) if its cohomological support loci are empty for every $i>0$and for every representation $x=\frac{a}{b}$. 
\end{defi}
Clearly we have that I.T.(0), implies M-regular which implies GV. In \cite{Caucci} and \cite{Ito22} various property of GV, M-regular and I.T.(0) sheaves are proven. Of most interest for us will be how these property are preserved under taking tensor products:
\begin{prop}[{\cite[Proposition 3.1, Theorem 3.2]{PP11},\cite[Prop. 3.4]{Caucci}}]\label{caucci}
Let $(X,l)$ a polarized abelian variety and  $\mathcal{F}\langle sl\rangle$ and $\mathcal{G}\langle tl\rangle$ two $\mathbb{Q}$-twisted sheaves such that one of them is locally free.
\begin{enumerate}
\item If $\mathcal{F}\langle sl\rangle$ is I.T.(0) and $\mathcal{G}\langle tl\rangle$ is GV, then $\mathcal{F}\otimes\mathcal{G}\langle(s+t)l\rangle$ is I.T.(0).
\item If $\mathcal{F}\langle sl\rangle$ and $\mathcal{G}\langle tl\rangle$ are both M-regular, then $\mathcal{F}\otimes\mathcal{G}\langle(s+t)l\rangle$ is M-regular.
\end{enumerate}
\end{prop}

\subsection{Cohomological rank functions on abelian varieties}\label{crf}

Let $X$ be an abelian variety, and $l\in\operatorname{NS}(X)$ a polarization. 
For an object $\mathcal{F}\in D^b(X)$  we can consider
$$
h^i_{\text{gen}}(X,\mathcal{F}),
$$
the dimension of the i-th (hyper)cohomology group $H^i(X,\mathcal{F}\otimes P_  \alpha)$ for the general $  \alpha\in\hat{X}$. Then the \emph{cohomological rank functions associated to $\mathcal{F}$ and $l$} is 
$$h^i_{\mathcal{F},l}:\mathbb{Q}\rightarrow\mathbb{Q}_{\ge 0}$$ defined via the assignment
$$\frac{a}{b}\mapsto b^{-2g}h^i_{\text{gen}}(X,b_X^*\mathcal{F}\otimes L^{ab}),$$
with $L$ any line bundle such that $c_1(L)=l$.\par
We see that the value of $h^i_{\mathcal{F},l}(x)$ depends only of $\mathcal{F}\langle xl\rangle$ (\cite[Remark 2.2]{JiPa2020}). Furthermore, if a $\mathbb{Q}$-twisted sheaf $\mathcal{F}\langle xl\rangle$ satisfies I.T.(0)  then we have that $h^i_{\mathcal{F},l}(x)=0$. In fact, let $x=\frac{a}{b}$ and chose $L$, such that $c_1(L)=l$. Then for all $\alpha\in\hat{X}$ we have that
$$
h^i(X,b_X^*\mathcal{F}\otimes L^{ab}\otimes\alpha)=0,
$$
when $i>0$.
In particular $h^i(X,b_X^*\mathcal{F}\otimes L^{ab})_{\text{gen}}=0$ and $h^i_{\mathcal{F},l}(x)=0$.
In addition to this, cohomological rank functions have been used to define invariants associated to polarizations of abelian varieties. Let $p\in X$ a closed point, we denote by $\mathcal{I}_p$ its associated sheaf of ideals. Given $l$ a polarization on $X$, we have that the quantity
\[\beta(l)\coloneqq\operatorname{inf}\{x\in\mathbb{Q}\:|\: h^1_{\mathcal{I}_p,l}(x)=0\}\]
does not depends from the choice of $p$. Following the work of Jiang--Pareschi and Caucci, we call it the \emph{basepoint-freenes treshold} of $l$. The reason behind the name is that $\beta(l)<1$ if and only if all representative $L$ of $l$
 are globally generated, thus this quantity, in some ways , measures how much representative of $l$ are globally generated.\par
 When $l$ is globally generated, then, for any line bundle $L$ with $c_1(L)=l$, we can consider the short exact sequence
 \begin{equation}\label{sesimportant}
     0\rightarrow M_L\longrightarrow H^0(X,L)\otimes\mathcal{O}_X\longrightarrow L\rightarrow 0
 \end{equation}
 as in \ref{sub:N_P}. We then define the invariant
 \[
 k(l)\coloneqq \inf\{x\in\mathbb{Q}^+\:|\: h^1_{M_L,l}(x)=0\}.
 \]
 When $L$ is very ample, this invariant measures the complexity of syzygies of the embedding induced by $L$, as shown by Caucci in \cite{Caucci}.\par
 By the work of Jiang--Pareschi \cite[Thm. D]{JiPa2020}, we have that the invariants $\beta(l)$ and $k(l)$ are related by the following equation
 \begin{equation}\label{eq:kappabeta}
     k(l)=\frac{\beta(l)}{1-\beta(l)}
 \end{equation}
  when $l$ is base-point free (that is when $\beta(l)<1$).
\section{Ampleness of Tautological Bundles}\label{sec3}

Let $Z$ be a smooth projective complex variety. Let $V$ be a vector bundle of rank $\mathfrak{r}$ on $Z$. Let $\mathbb{P}_Z(V)$ be the associated projective bundle with the projection morphism $\pi:\mathbb{P}_Z(V)\rightarrow Z$. Let $\mathcal{O}_{\mathbb{P}(V)}(1)$ be the tautological line bundle of $\mathbb{P}_Z(V)$.
\begin{defi} 
A vector bundle $V$ on $Z$ is \emph{ample} (respectively \emph{very ample} or \emph{nef}) if the tautological bundle $\mathcal{O}_{\mathbb{P}(V)}(1)$ is ample (respectively very ample or nef) on the projective bundle $\mathbb{P}_Z(V)$. 
\end{defi}
A generalization of this notion is the concept of $k$-jet ampleness \cite[\S 1.2]{BDS}.

\begin{defi}[$k$-jet ampleness]
Let $k\geq 0$. The bundle $V$ is $k$-jet ample if the evaluation map 
\[
H^0(Z,V)\rightarrow H^0(Z,V\otimes (\mathcal{O}_Z/(m_{x_1}^{k_1}\otimes \cdots\otimes m_{x_p}^{k_p})))=\bigotimes_{i=1}^pH^0(Z,(\mathcal{O}_Z/m_{x_i}^{k_i}))
\]
is surjective for any choice of $t$ distinct points $x_1,\ldots,x_p\in Z$ such that $\sum_{i=1}^pk_i=k+1$, where $m_{x_i}$ is the maximal ideal sheaf of $x_i\in Z$.
\end{defi}

A vector bundle $V$ is $0$-jet ample if and only if it is globally generated. In addition, according to \cite[Proposition 4.2]{BDS}, if $V$ is a $1$-jet vector bundle on $Z$, then $V$ is very ample on $Z$.

Given a very ample vector bundle $V$ on $Z$, it make sense to give conditions ensuring property $N_p^r$ holds for the embedding induced by $\mathcal{O}_{\mathbb{P}(V)}(1)$. In \cite{Park} a cohomological criterion for property $N_p^r$ for the tautological bundle of a projective bundle $\mathbb{P}(V)\rightarrow Z$ is given. Before recalling the statement we need to fix some notation.\par
To this aim let $V$ be a very ample vector bundle on $Z$. In particular, we know that $V$ is globally generated and we have a short exact sequence
\[ 0\rightarrow M_{V}\longrightarrow H^0(Z, V)\otimes\mathcal{O}_Z\longrightarrow V\rightarrow 0,\]
where $M_V$ denotes the kernel of the evalutaion map $H^0(Z,V)\otimes\mathcal{O}_Z\rightarrow V$.

\begin{prop}\label{PARK}
    The tautological line bundle $\mathcal{O}_{\mathbb{P}(V)}(1)$ satisfies property $N_p^r$
 if 
 \[H^k\left(Z,\bigwedge^iM_V\otimes V^{\otimes h}\right)=0\]
 for every $1\leq i\leq p+1 $, $k\geq 1$ and $h\geq r+1.$
 \end{prop}
 \begin{proof}
    This is essentially \cite[Prop. 3.2]{Park}. Let $X=\mathbb{P}(V)$ and denote by $L$ the tautologiacal bundle on $X$, and let $\pi:X\rightarrow Z$ the morphism giving the projective bundle structure. It is well known that property $N_p^r$ is implied by the vanishing 
      \[H^1\left(X,\bigwedge^i M_L\otimes V^{\otimes h}\right)=0\]
      for every $h\geq r+1$ and for every $1\leq i\leq p+1$. In the proof of \cite[Prop. 3.2]{Park}, it is shown how that this is implied by the vanishing
     \begin{equation}\label{parkexplained}H^1\left(X,\bigwedge^l\pi^*M_V\otimes\Omega^m_{X/Y}\otimes L^{\otimes h}\right)=0,
     \end{equation}
      for every $l+m=i$ and for $h\geq r+1$.
      As a consequence of projection formula plus an easy cohomological argument, we have that \eqref{parkexplained}  is implied  by the vanishing 
      \begin{equation}\label{parkexplained2}
H^k\left(Z,\bigwedge^lM_V\otimes\bigwedge^{m+1+k}V\otimes S^{m+h-k} V\right)=0,
      \end{equation}
       for every $k\geq 1$, $l+m=i$, and for $h\geq r+1$. Now, 
      $\bigwedge^{m+1+k}V\otimes S^{m+j-k}V$ is a summand of $V^{\otimes m+j}$, therefore we see that \eqref{parkexplained2}
 is implied      \[H^k\left(Z,\bigwedge^lM_V\otimes V^{\otimes m+h}\right)=0\]
      for every $l+m=i$ and every $h\geq r+1$.
 \end{proof}
\section{Projective Bundles on Abelian Varieties}\label{sec4}
Before proceeding with the proof of our main results we recall our specific setting and introduce some notation. Let $(X,l)$ be a complex abelian variety and consider a vector bundle $E$ of rank $\mathfrak{r}$ on $X$ admitting a filtration
\begin{equation}\label{eq:filtration}
    0=E_0\subset E_1\subset E_2\subset\cdots\subset E_\mathfrak{r}=E,
\end{equation}
such that the quotients
$\beta_t:=E_t/E_{t-1}$ lie in $\Pic^0(X)$, for every $t=1,\ldots,\mathfrak{r}$.
Let $Y$ be the projective bundle $\mathbb{P}_X(E)$ with associated projective morphism $\pi:Y\rightarrow X$. Given a line bundle $A$ with $\mathrm{c}_1(A)=\ell$ denote by $\pi_A:Y_A\rightarrow X$ be the projective bundle associated to the vector bundle $E\otimes A$. Recall that there is an isomorphism 
\[\varphi:Y \rightarrow Y_A\]
such that $\varphi^*\mathcal{O}_{Y_A}(1)\simeq \mathcal{O}_{Y}(1)\otimes \pi^*A$, thus to Theorem \ref{TheoremA} it will be enough to show that $\mathcal{O}_{Y_A}(1)$ is very ample and satisfies property $N_p^r$. 
\subsection{Very ampleness and vanishing properties of $E\otimes A$}
By definition the very ampleness of $\mathcal{O}_{Y_A}(1)$ is equivalent to the very ampleness of the vector bundle $E\otimes A$ on $X$, thus the very ampleness part of our Theorem \ref{TheoremA} is a direct consequence of the following.
\begin{lemma}\label{lem:kample}
    Suppose that $X$ and $E$ are as in the statement of the Theorem \ref{TheoremA}. Let $\ell$ be a polarization on $X$ satisfying $\beta(\ell)<\frac{1}{p+2}$. Then for every $A$ ample line bundle with $\mathrm{c}_1(A)=\ell$ we have that the vector bundles $E_t\otimes A$ are $(p+1)$-jet ample. In particular $E\otimes A$ is $(p+1)$-jet ample.
\end{lemma}
\begin{proof}
    By \cite[Theorem D]{Caucci} we have that $A\otimes \beta_t$ is $(p+1)$-jet ample for every $t=1,\cdots, \mathfrak{r}$. In particular $E_1\otimes A$ is $(p+1)$-jet ample. By using the short exact sequence
    \begin{equation}\label{eqn4.2}
    0\rightarrow E_{t-1}\otimes A\longrightarrow E_t\otimes A\longrightarrow \beta_t\otimes A\rightarrow 0
    \end{equation}
    and the fact that extension of $(p+1)$-jet ample are $(p+1)$-jet ample \cite[Proposition 2.6]{BDS}, we can show inductively that $E_t\otimes A$ is $(p+1)$-jet ample.
\end{proof}
With a very similar argument to the one above, one can also show the following useful result.
\begin{lemma}\label{lem:EAIT0}
    Suppose that $X$ and $E$ are as in the statement of the Theorem. Then for every $A$ ample line bundle with on $X$ we have that the vector bundles $E_t\otimes A$ satisfy I.T.(0). In particular $E\otimes A$ satisfies I.T.(0).
\end{lemma}
\begin{proof}
By the short exact sequence \eqref{eqn4.2} and the fact that extension of I.T.(0) sheaves, we can show inductively that $E\otimes A$ is I.T.(0).
\end{proof}

\subsection{Syzygies for tautological bundles} In this paragraph we want to show that $\mathcal{O}_{Y_A}(1)$ satisfies property $N_p^r$ when $\beta(\ell)<\min\left\{\frac{r+1}{p+r+2},\frac{1}{2}\right\}$. As a consequence of Lemma \ref{lem:kample} we have that $E\otimes A$ is globally generated whenever $\beta(\mathrm{c}_1(A))<\frac{1}{2}$. In particular,  we have a short exact sequence
\begin{equation}\label{eq:sesEA}
    0\rightarrow M_{E\otimes A}\longrightarrow H^0(X, E\otimes A)\longrightarrow E\otimes A\rightarrow 0
\end{equation}
By Proposition \ref{PARK} , to show that property $N_p^r$ holds for $\mathcal{O}_{Y_A}(1)$ we need to show the vanishing
\begin{equation}\label{eq:vanishing 1}
    H^k\left(X,\bigwedge^jM_{E\otimes A}\otimes (E\otimes A)^{\otimes h}\right)=0
\end{equation}
for every $k\geq 1$, $h\geq r+1$ and every $1\leq j\leq p+1$.
Since we are working in characteristic 0 we have that the vanishing \eqref{eq:vanishing 1} is implied by the following one:
\begin{equation}\label{eq:vanishing 1}
    H^k(X,M_{E\otimes A}^{\otimes j}\otimes (E\otimes A)^{\otimes h})=0
\end{equation}
which in particular is verified when the vector bundle $M_{E\otimes A}^{\otimes j}\otimes (E\otimes A)^{\otimes h}$ satisfies I.T.(0). Thus we aim to showing the following result. We note that the second statement (2) from the proposition below is essential to prove the first statement (1). 
\begin{prop}\label{prop:ITO}
Suppose that $X$ and $E$ are as in the statement of the Theorem. Let $\ell$ be a polarization on $X$ satisfying 
\[
\beta(\ell)<\min\left\{\frac{r+1}{p+r+2},\frac{1}{2}\right\},
\]
then, for every $A$ ample line bundle with $\mathrm{c}_1(A)=\ell$ we have that 
\begin{enumerate}
\item the vector bundle $M_{E\otimes A}^{\otimes j}\otimes (E\otimes A)^{\otimes h}$ satisfies I.T.(0) for every $h\geq r+1$ and for every $1\leq j\leq p+1$, and
\item the vector bundle 
       \begin{equation}\label{eqn1rd}
    \bigotimes_{t=1}^{\mathfrak{r}}M_{\beta_t\otimes A}^{\otimes k_t}\otimes E^{\otimes u}\otimes \beta_{\mathfrak{r}}^{\otimes h-u}\otimes A^{\otimes h} \end{equation}
       satisfies I.T.(0) whenever $\sum k_t\leq p+1$ and $h\geq r+1$ for $0\leq u\leq h$.
\end{enumerate}
\end{prop}
The proof of this result will take the reminder of this section. We start with the following lemma for the case $u=0$ of (2) in Proposition \ref{prop:ITO}.

\begin{lemma}\label{IT0Q}
    Let $(X,\ell)$ be a complex polarized abelian variety with 
    \[
    \beta(\ell)<\min\left\{\frac{r+1}{p+r+2},\frac{1}{2}\right\}.
    \]
    Let $A$ an ample line bundle representing the polarization $\ell$, and take $k_1,\cdots, k_\mathfrak{r}$  non negative integers such that  $\sum_tk_t\leq p+1$ and $\beta_t$ with $t=1,\ldots, \mathfrak{r}$ element in $\Pic^0(X)$. Then 
    \[M_{\beta_1\otimes A}^{\otimes k_1}\otimes M_{\beta_2\otimes A}^{\otimes k_2}\otimes\cdots\otimes M_{\beta_{\mathfrak{r}}\otimes A}^{\otimes k_\mathfrak{r}}\langle h\ell\rangle \]
    satisfies I.T.(0) for all $h\geq r+1$.
\end{lemma}
\begin{proof}
   Since $\beta_t\otimes A$ has class equal to $\ell$, our assumptions warrant that $k(\ell)<\frac{r+1}{p+1}$, and in particular $M_{\beta_t\otimes A}\left\langle\frac{r+1}{p+1}\ell\right\rangle$ satisfy I.T.(0). By Proposition \ref{caucci}, powers $\left(M_{\beta_t\otimes A}\left\langle\frac{r+1}{p+1}\ell\right\rangle\right)^{\otimes k_t}$ satisfy I.T.(0). Now we write
   \begin{align*}
       M_{\beta_1\otimes A}^{\otimes k_1}\otimes \cdots\otimes M_{\beta_{\mathfrak{r}}\otimes A}^{\otimes k_\mathfrak{r}}\langle h\ell\rangle
       &= M_{\beta_1\otimes A}^{\otimes k_1}\otimes \cdots\otimes M_{\beta_{\mathfrak{r}}\otimes A}^{\otimes k_\mathfrak{r}}\left\langle \left(\frac{(r+1)\cdot\sum k_t}{p+1}+h-\frac{(r+1)\cdot\sum k_t}{p+1}\right)\ell\right\rangle\\
       &= \bigotimes_{t=1}^\mathfrak{r}M_{\beta_t\otimes A}^{\otimes k_t}\left\langle\frac{(r+1)\cdot k_t}{p+1}\ell\right\rangle\otimes \mathcal{O}_X\left\langle \left(h-\frac{(r+1)\cdot\sum k_t}{p+1}\right)\ell\right\rangle\\
       &=\bigotimes_{t=1}^\mathfrak{r} \left(M_{\beta_t\otimes A}\left\langle\frac{r+1}{p+1}\ell\right\rangle\right)^{\otimes k_t}\otimes \mathcal{O}_X\left\langle \left(h-\frac{(r+1)\cdot\sum k_t}{p+1}\right)\ell\right\rangle.
   \end{align*}
   Since $\sum k_t\leq p+1$, we have that $h-\dfrac{(r+1)\cdot\sum k_t}{p+1}\geq 0$ whenever $h\geq r+1$. In particular the $\mathbb{Q}$-twisted sheaf $\mathcal{O}_X\left\langle \left(h-\frac{(r+1)\cdot\sum k_t}{p+1}\right)\ell\right\rangle$ satisfies GV. On the other side $\bigotimes_{t=1}^{\mathfrak{r}} \left(M_{\beta_t\otimes A}\left\langle\frac{r+1}{p+1}\ell\right\rangle\right)^{\otimes k_t}$ is a product of I.T.(0) $\mathbb{Q}$-twisted vector bundles, and therefore it satisfies I.T.(0). We conclude by applying Proposition \ref{caucci}.\end{proof}

Let us proceed with the proof of the second statement in Proposition \ref{prop:ITO}. 
   \begin{proof}[Proof of Proposition \ref{prop:ITO} (2)]  
  We prove this by induction on $u$. Let $\alpha\in \Pic^0(X)$ and set $u=0$. 
  By Lemma \ref{IT0Q}, we have that 
\[H^k\left(X,\bigotimes_{t=1}^\mathfrak{r} M_{\beta_t\otimes A}^{\otimes k_t}\otimes (\beta_{\mathfrak{r}}\otimes A)^{\otimes h}\otimes\alpha\right)=0\]
for every $k\geq 1$. We deduce immediately that 
$\bigotimes_{t=1}^\mathfrak{r} M_{\beta_t\otimes A}^{\otimes k_t}\otimes (\beta_{\mathfrak{r}}\otimes A)^{\otimes h}$ satisfies I.T.(0). 

We now assume that 
      \begin{equation}\label{eqn4.5}\bigotimes_{t=1}^{\mathfrak{r}}M_{\beta_t\otimes A}^{\otimes k_t}\otimes E^{\otimes u-1}\otimes \beta_{\mathrm{r}}^{\otimes h-u+1}\otimes A^{\otimes h} \end{equation}
      satisfies I.T.(0) for some $u\geq 1$.    
We will show the more general statement that all the vector bundles 
\begin{equation}\label{eqn:2nd}
\bigotimes_{t=1}^{\mathfrak{r}}M_{\beta_t\otimes A}^{\otimes k_t}\otimes E^{\otimes u-1}\otimes E_{s}\otimes \beta_{\mathrm{r}}^{\otimes h-u}\otimes A^{\otimes h}
\end{equation}
 satisfies I.T.(0) for any $ s\geq 1$, and we prove this by induction on $s$. 
 
 Let $s=1$. After taking the arbitrary $\alpha$ to be $E_1^\vee\otimes \beta_{\mathfrak{r}}\otimes \alpha\in\Pic^0(X)$ and by the induction hypothesis \eqref{eqn4.5} from the first claim, we obtain that 
 \begin{align*}
 H^k\Bigg(X,\bigotimes_{t=1}^\mathfrak{r} M_{\beta_t\otimes A}^{\otimes k_t}\otimes & E^{\otimes u-1}\otimes E_{1}\otimes \beta_{\mathrm{r}}^{\otimes h-u}\otimes A^{\otimes h}
\otimes\alpha\Bigg)\\
&=H^k\left(X,\bigotimes_{t=1}^\mathfrak{r} M_{\beta_t\otimes A}^{\otimes k_t}\otimes E^{\otimes u-1}\otimes \beta_{\mathrm{r}}^{\otimes h-u+1}\otimes A^{\otimes h}
\otimes\alpha\right)=0\end{align*}
for all $k$. This immediately gives us that $M_{\beta_t\otimes A}^{\otimes k_t}\otimes  E^{\otimes u-1}\otimes E_{1}\otimes \beta_{\mathrm{r}}^{\otimes h-u}\otimes A^{\otimes h}$ satisfies I.T.(0).

Let us assume that 
\begin{equation}\label{eqn4.6}
\bigotimes_{t=1}^{\mathfrak{r}}M_{\beta_t\otimes A}^{\otimes k_t}\otimes E^{\otimes u-1}\otimes E_{s-1}\otimes \beta_{\mathrm{r}}^{\otimes h-u}\otimes A^{\otimes h}
\end{equation}
satisfies I.T.(0) for some $s\geq 2$. We look at the short exact sequence
\begin{align*}
0\rightarrow \bigotimes_{t=1}^{\mathfrak{r}}M_{\beta_t\otimes A}^{\otimes k_t}\otimes E^{\otimes u-1}\otimes E_{s-1}\otimes \beta_{\mathrm{r}}^{\otimes h-u}\otimes A^{\otimes h}\longrightarrow \bigotimes_{t=1}^{\mathfrak{r}}M_{\beta_t\otimes A}^{\otimes k_t}\otimes E^{\otimes u}\otimes E_s\otimes\beta_{\mathrm{r}}^{\otimes h-u}\otimes A^{\otimes h}\\
\longrightarrow \bigotimes_{t=1}^{\mathfrak{r}}M_{\beta_t\otimes A}^{\otimes k_t}\otimes E^{\otimes u-1}\otimes \beta_{\mathrm{r}}^{\otimes h-u+1}\otimes A^{\otimes h}\rightarrow 0.
\end{align*}
As an extension of I.T.(0) sheaves, we can conclude that $\bigotimes_{t=1}^{\mathfrak{r}}M_{\beta_t\otimes A}^{\otimes k_t}\otimes E^{\otimes u}\otimes E_s\otimes\beta_{\mathrm{r}}^{\otimes h-u}\otimes A^{\otimes h}$ satisfies I.T.(0) by the induction hypothesis \eqref{eqn4.5} and \eqref{eqn4.6}. By taking $s=\mathfrak{r}$ from \eqref{eqn:2nd}, we have the statement \eqref{eqn1rd}.
  \end{proof}

We finish this section with the proof of Proposition \ref{prop:ITO} (1).

   \begin{proof}[Proof of Proposition \ref{prop:ITO} (1)]  
   After taking $u=h$ from \eqref{eqn1rd} from Proposition \ref{prop:ITO} (2), we have that 
   \[H^k\left(X,\bigotimes_{t=1}^\mathfrak{r} M_{\beta_t\otimes A}^{\otimes k_t}\otimes (E\otimes A)^{\otimes h} \otimes\alpha\right)=0\] for every $\alpha\in\Pic^0(X)$ and every $k\ge 1$, and every non-negative integers $k_t$ such that $\sum k_t\leq p+1$. Let $j=\sum k_t$. By \cite[Lemma 4.4]{Chi19}, we have that 
   \[H^k\left(X,M_{E\otimes A}^{\otimes j}\otimes (E\otimes A)^{\otimes h} \otimes\alpha\right)=0\] for every $\alpha\in\Pic^0(X)$, every $k\ge 1$, $h\geq r+1$, and $1\leq j\leq p+1$ as we wanted.
   \end{proof}

     We observe that Corollary \ref{introcor1} follows immediately from Theorem \ref{TheoremA} with $r=0$.

   \begin{cor}\label{introcor1}
    Let $X$ be an abelian varieties, $E$ a vector bundle as in the statement of Theorem \ref{TheoremA}, and $A$ an ample line bundle on $X$, then $\mathcal{O}_Y(1)\otimes\pi^* A^{\otimes n}$ satisfies $N_p$ for every $n\geq p+3$.
  \end{cor}
  
  In fact, if $A$ is an ample line bundle, we have that 
   $$\beta(n\operatorname{c}_1(A))=\frac{1}{n}\beta(\operatorname{c}_1(A))\leq \frac{1}{n}.$$
In particular, if $n\geq p+3$ we have that $\beta(n\operatorname{c}_1(A)))<\frac{1}{p+2}$ and therefore $\mathcal{O}_{\mathbb{P}(E)}(1)\otimes\pi^*A^{\otimes n}$ satisfies $N_p$.

\subsection{Finer results for higher syzygies using the work of Ito}
Let $\mathcal{I}_o\subseteq\mathcal{O}_X$ be the maximal ideal sheaf corresponding to the origin $o$ in an abelian variety $X$. 
Recall that $A$ is an ample line bundle on $X$ with its polarization $\ell=c_1(A)$.

\begin{prop}\label{prop5.1}
Suppose $(X,E)$ is the pair of abelian variety $X$ and vector bundle $E$ of rank $\mathfrak{r}$. Let $p\geq 0$.
If
$
\mathcal{I}_o\left\langle \frac{1}{p+2}\ell\right\rangle
$ is M-regular,
 then $E\otimes A$ is a $(p+1)$-jet ample vector bundle on $X$. 
\end{prop}

\begin{proof}
By \cite[Theorem 1.6]{Ito22}, $\beta_t\otimes A$ and $E_1\otimes A$ are $(p+1)$-jet ample on $X$. The remaining arguments follow the same line of reasoning as the one used for Lemma \ref{lem:kample}, leading to its conclusion.
\end{proof}
As a corollary, we have that $E\otimes A$ is globally generated if $
\mathcal{I}_o\left\langle \frac{1}{2}\ell\right\rangle
$ is M-regular.

Now, we investigate the $N_p$-property of projective bundles over an abelian variety $X$. 

\begin{thm}\label{thm:5.4}
Suppose $(Y,X,\pi)$ as before. 
Let $A$ be an ample line bundle on $X$. Let $p\geq 1$.
 If $\mathcal{I}_o\left\langle \frac{r+1}{p+r+2}\ell\right\rangle$ is M-regular, then $\mathcal{O}_Y(1)\otimes\pi^*A$ satisfies $N_p^r$-property.  
 In particular, if $\mathcal{I}_o\left\langle \frac{1}{p+2}\ell\right\rangle$ is M-regular, then $\mathcal{O}_Y(1)\otimes\pi^*A$ satisfies $N_p$ property. 
\end{thm}
\begin{proof}
 We consider ample line bundles $\beta_t\otimes A$ for $1\leq t\leq \mathfrak{r}$. Since $\mathcal{I}_o\left\langle\frac{r+1}{p+r+2}\ell\right\rangle$ is M-regular, $M_{\beta_t\otimes A}\left\langle\frac{r+1}{p+1}\ell\right\rangle$ is also M-regular for $1\leq t\leq\mathfrak{r}$ by \cite[Proposition 4.1]{Ito22}. 
 We take $k_1$ to be an integer greater than or equal to 2 such that $2\leq k_1\leq p+1$. Then by \cite[Proposition 6.2]{Ito22}, $M_{\beta_1\otimes A}^{\otimes k_1}\left\langle\frac{(r+1)\cdot k_1}{p+1}\ell\right\rangle$ is I.T.(0). We take any nonnegative integers $k_2,\ldots,k_{\mathfrak{r}}$ such that $\sum_tk_t\leq p+1$. By Proposition \ref{caucci}
 $\bigotimes_{t=2}^\mathfrak{r} \left(M_{\beta_t\otimes A}\left\langle\frac{r+1}{p+1}\ell\right\rangle\right)^{\otimes k_t}$ is M-regular, and thus GV. Recall that $\mathcal{O}_X\left\langle \left(h-\frac{(r+1)\cdot\sum k_t}{p+1}\right)\ell\right\rangle$ is GV for $h\geq r+1$ from the proof of Proposition \ref{IT0Q}. By Proposition \ref{caucci}, we have that    
$M_{\beta_1\otimes A}^{\otimes k_1}\otimes M_{\beta_2\otimes A}^{\otimes k_2}\otimes\cdots\otimes M_{\beta_{\mathfrak{r}}\otimes A}^{\otimes k_\mathfrak{r}}\langle h\ell\rangle=\bigotimes_{t=1}^\mathfrak{r} \left(M_{\beta_t\otimes A}\left\langle\frac{r+1}{p+1}\ell\right\rangle\right)^{\otimes k_t}\otimes \mathcal{O}_X\left\langle \left(h-\frac{(r+1)\cdot\sum k_t}{p+1}\right)\ell\right\rangle$ satisfies I.T.(0) for all $h\geq r+1$. From the same arguments for Proposition \ref{prop:ITO}, we have that the vector bundle $M_{E\otimes A}^{\otimes j}\otimes (E\otimes A)^{\otimes h}$ satisfies I.T.(0) for all $h\geq r+1$ and for any $2\leq j\leq p+1$, where $j=\sum_tk_t$.

Now we consider the case when $j=1$. We set $k_s=1$ and the other $k_i=0$ for every $i$ except $s$, (i.e., $i\neq s$). 
Let us look at the following decompositions
   \begin{align*}
       M_{\beta_s\otimes A}\left\langle h\ell\right\rangle&
     =  M_{\beta_s\otimes A}\left\langle ((r+1)+ h-(r+1))\ell\right\rangle
       = M_{\beta_s\otimes A}\langle (r+1)\ell\rangle\otimes\mathcal{O}_X\left\langle \left( h-r-1)\right)\ell\right\rangle\\
       &= M_{\beta_s\otimes A}\left\langle\frac{r+1}{p+1}(p+1)\ell\right\rangle\otimes \mathcal{O}_X\left\langle \left(h-r-1)\right)\ell\right\rangle.
   \end{align*}
   Since $p+1\geq 2$, by \cite[Proposition 6.2]{Ito22} we have that $M_{\beta_s\otimes A}\left\langle\frac{r+1}{p+1}(p+1)\ell\right\rangle$ is I.T.(0). In addition, the $\mathbb{Q}$-twisted sheaf $\mathcal{O}_X\left\langle \left(h-r-1)\right)\ell\right\rangle$ satisfies GV, as $h-r-1\geq 0$. Thus, the product $M_{\beta_s\otimes A}\left\langle h\ell\right\rangle$ satisfies I.T.(0). 
   
   Combining all together, we get that $M_{E\otimes A}^{\otimes j}\otimes (E\otimes A)^{\otimes h}$ satisfies I.T.(0) for all $h\geq r+1$ and for every $1\leq j\leq p+1$.

We can deduce the last statement as the special case when $r=0$. 
\end{proof}

The theorem above provides an analogous result to the case of ample line bundles in \cite[Theorem 1.5]{Ito22}. 
The following corollary is a generalization of the consequences presented in Theorem \ref{thm:5.4}.
\begin{cor}
Suppose $(Y,X,\pi)$ as before. 
Let $A$ be an ample line bundle on $X$. Let $p\geq 1$.
\begin{enumerate}
\item If $\mathcal{I}_o\left\langle \frac{n}{p+2}\ell\right\rangle$ is M-regular, then $\mathcal{O}_Y(1)\otimes\pi^*A^{\otimes n}$ satisfies $N_p$-property. 
\item If $\mathcal{I}_o\left\langle \frac{n(r+1)}{p+r+2}\ell\right\rangle$ is M-regular, then $\mathcal{O}_Y(1)\otimes\pi^*A^{\otimes n}$ satisfies $N_p^r$-property. 
\end{enumerate}
\end{cor}

We also have the following corollary regarding the condition on $n$ for projective bundles over abelian varieties to have $N_p$ and $N_p^r$ properties. 

\begin{cor}\label{cor:5.5}
Let $p\geq 1$ be an integer. Let $A$ be an ample line bundle on $X$ with no base divisor. 
\begin{enumerate}
\item if $n\geq p+2$, then $\mathcal{O}_Y(1)\otimes\pi^*A^{\otimes n}$ satisfies $N_p$-property, and
\item More generally, if $n\geq (p+r+2)/(r+1)$, then $\mathcal{O}_Y(1)\otimes\pi^*A^{\otimes n}$ satisfies $N_p^r$-property for $r\geq 0$. 
\end{enumerate}
\end{cor}
\begin{proof}
(1) According to \cite[Remark 3.6]{PP04}, $L$ does not have a base divisor if and only if $\mathcal{I}_o\langle l\rangle$ is M-regular for $c_1(L)=l$. So, given $L=E_t/E_{t-1}\otimes A^{\otimes n}$ for $1\leq t\leq \mathfrak{r}$, if $A$ does not have a base divisor, then $\mathcal{I}_o\left\langle\frac{1}{p+2}l\right\rangle=\mathcal{I}_o\left\langle\frac{n}{p+2}\ell\right\rangle$ is M-regular if ${n}/({p+2})\geq 1$, where $c_1(A)=\ell$. Thus, we establish (1) using the same argument as in Theorem \ref{thm:5.4}. Similarly, the second statement (2) holds.
\end{proof}

We note that Corollary \ref{cor:5.5} (1) generalizes Corollary \ref{introcor1} when $p\geq 1$.

\section{Application to projective bundles on desinglarization of compactified Jacobians}\label{sec5}

In this section, we investigate syzygies of projective bundles on normalization of compactified Jacobians. Our application provides an improvement of the result presented in \cite[Theorem 5.4]{Chi19} when it comes to an ample line bundle with no base divisor. 

We consider an irreducible nodal curve $X$ over a complex number $\mathbb{C}$ with a node $y$. Let $N$ be the desingulazation (normalization) of $X$, with the normalization map $p:X_0\rightarrow X$. We denote by $x$ and $z$ the inverse image of the node $y$ so that $p^{-1}(y)=\{x,z\}$. If the genus of $X_0$ is $g$, then the arithmetic genus $p_a(X)$ of $X$ is $g+1$. 

Let $d\in \mathbb{Z}$ and $J^d(X)$ the Jacobian of line bundles, locally free sheaves of rank $1$, of degree $d$ on $X$. The Jacobian $J^d(X)$ admits a natural compactification $\overline{J^d}(X)$ which parametrizes torsion-free rank 1 sheaves of degree $d$ on $X$. Let $h:\widetilde{J^d}(X)\rightarrow\overline{J^d}(X)$ be a natural desingularization of $\overline{J^d}(X)$. 
In particular, $\widetilde{J^d}(X)$ can be described as a projective bundle as follows.

We let $\mathcal{P}$ be the Poincar\'{e} bundle on $J^d(X_0)\times X_0$ such that for some $c\in X_0$,
$\mathcal{P}|_{J^d(X_0)\times\{c\}}$ is trivial. We denote by $\mathcal{P}_{x}$ and $\mathcal{P}_{z}$ the trivial line bundles from the restriction of $\mathcal{P}$ on $J^d(X_0)\times\{x\}$ and $J^d(X_0)\times\{z\}$ respectively. Then we have
\[
\widetilde{J^d}(X)\cong\mathbb{P}(\mathcal{P}_x\oplus\mathcal{P}_z)
\]
by \cite[Proposition 2.1]{BP08}. So, $\widetilde{J^d}(X)$ can be viewed as a $\mathbb{P}^1$-bundle over $J^d(X_0)$. For more details about the construction of $\widetilde{J^d}(X)$, readers may refer to \cite{BP08} and \cite[Proposition 12.1]{OS79}.
Let $E$ be the direct sum $\mathcal{P}_x\oplus\mathcal{P}_z$ of degree $0$ line bundles. 

%

Now, we are ready to have the following theorem. 
\begin{thm}
Let $X$ and $E$ be as above, and $A$ an ample line bundle with ${c}_1(A)=\ell$ on $J(X_0)$ with no base divisor. Let $\pi:\widetilde{J}(X)\rightarrow J(X_0)$ be the projection morphism with $\mathcal{O}_{\mathbb{P}(E)}(1)$ the tautological line bundle.  Let $p\geq 1$.

\begin{enumerate}
\item If $n\geq p+2$, then $\mathcal{O}_{\mathbb{P}(E)}(1)\otimes\pi^*A^{\otimes n}$ satisfies $N_p$-property, and
\item More generally, if $n\geq (p+r+2)/(r+1)$, then $\mathcal{O}_{\mathbb{P}(E)}(1)\otimes\pi^*A^{\otimes n}$ satisfies $N_p^r$-property for $r\geq 0$. \end{enumerate}
\end{thm}
\begin{proof}
Since $E=\mathcal{P}_x\oplus\mathcal{P}_z$ is a rank $2$ vector bundle admitting a filtration as in \eqref{eq:filtration} with $t=2$, we can apply Corollary \ref{cor:5.5} to have (1) and (2).
\end{proof}

%
%
%
%
\section*{Acknowledgment}
Both authors are grateful to the anonymous referee whose comments greatly improved the exposition. \par

This project has started during Virtual Workshop II for Women in Commutative Algebra and Algebraic Geometry hosted by the Fields Institute and organized by M. Harada and C. Miller. We are really grateful to the Fields Institute, C. Miller and M. Harada for this opportunity to make new connection and create collaborations. We, in addition, wish to thanks both the organizers and the participants to the workshop for creating a friendly and stimulating environment, and for the engaging mathematical discussion.

In addition both authors are grateful to A. Ito for his valuable comments to a previous version of this manuscript and for pointing out that our Theorem A was a particular case of a result of his. 

 M. Jeon is partially supported by AMS-Simons Travel Grant. S. Tirabassi was partially supported by the Knut and Alice Wallenberg Foundation under grant no. KAW 2019.0493 and the VR grant 2023-03837.
\bibliographystyle{alpha}
\bibliography{biblio}
\end{document}